\documentclass[10pt]{amsart}
\usepackage[psamsfonts]{amssymb}
\usepackage{amsthm,amscd}
\usepackage{type1cm}
\usepackage{color}

\newtheorem{thm}{Theorem}

\newtheorem{lem}{Lemma}
\newtheorem{cor}{Corollary}
\newtheorem{pro}{Proposition}
\theoremstyle{remark}

\theoremstyle{definition}

\allowdisplaybreaks[4]

\title[2-dimensional endo-commutative straight algebras of type I]{A classification of 2-dimensional endo-commutative straight  algebras of type I}

\author[S.-E. Takahasi]{Sin-Ei Takahasi}
\author[K. Shirayanagi]{Kiyoshi Shirayanagi}
\author[M. Tsukada]{Makoto Tsukada}
\address[S.-E. Takahasi]{Laboratory of Mathematics and Games\\ Katsushika 2-371\\ Funabashi\\ Chiba 273-0032\\ Japan}
\address[K. Shirayanagi, M. Tsukada]{Department of Information Science\\ Toho University\\ Miyama 2-2-1\\ Funabashi\\ Chiba 274-8510\\ Japan}
\email{sin\_ei1@yahoo.co.jp}
\email[K. Shirayanagi (Corresponding author)]{kiyoshi.shirayanagi@is.sci.toho-u.ac.jp}
\email[M. Tsukada]{tsukada@is.sci.toho-u.ac.jp}

\makeatletter
\@namedef{subjclassname@2020}{\textup{2020} Mathematics Subject Classification}
\makeatother


\begin{document}

\begin{abstract}
  In this paper, we provide a complete classification of 2-dimensional endo-commutative straight algebras of type I over any field.
  An endo-commutative algebra is a non-associative algebra in which the square mapping preserves multiplication.  Type I denotes a distinguishing characteristic of
  its structure matrix of rank 2.
  We list all multiplication tables of these algebras up to isomorphism.
\end{abstract}
\maketitle
\section{Introduction}\label{sec:intro}
Let $A$ be a non-associative algebra.  The square mapping $x\mapsto x^2$ from $A$ to itself yields various important concepts of $A$.  We define $A$ to be {\it{endo-commutative}}, if the square mapping of $A$ preserves multiplication, that is, $x^2y^2=(xy)^2$ holds for all $x, y\in A$.
This terminology stems from the identity $(xx)(yy)=(xy)(xy)$, which illustrates the property of {\it inner} commutativity.
The concept of endo-commutativity was first introduced by the authors \cite{TST1}, where we provided a complete classification of endo-commutative algebras of dimension 2 over the trivial field of two elements.

We categorize two-dimensional algebras into two distinct types: {\it{curled}} and {\it{straight}}. Specifically, a 2-dimensional algebra is curled if the square of any element $x$ is a scalar multiple of $x$, otherwise it is straight.  In our previous work \cite{TST2}, we presented a complete classification of 2-dimensional endo-commutative curled algebras over any non-trivial field.  Research related to curled algebras can be found in \cite{level2,length1}.  In addition, in \cite{TST3}, we provided a complete classification of 2-dimensional endo-commutative straight algebras of rank 1 over any non-trivial field, where the rank of an algebra is defined as the rank of its structure matrix, which represents the matrix of structure constants with respect to a linear base.  

Our objective is to achieve a complete classification of 2-dimensional endo-commutative straight algebras of rank 2 over any field. To accomplish this,
we categorize these algebras into three types: type I, type II, and type III
based on a distinguishing characteristic of their structure matrix.

We can find classifications of associative algebras of dimension 2 over the real and complex number fields in \cite{onkochishin}. For other studies on 2-dimensional algebras, see \cite{moduli,2dim,variety,classification}. 

To begin, we proceed with the characterization of 2-dimensional endo-commutative straight algebras of type I over any field.
This characterization involves the structure matrix, which is determined by the product of each pair of elements from a linear base. Proposition 1 provides comprehensive details on this characterization.
As a result of this characterization, we discover that such algebras are commutative, but neither unital nor associative.
Corollary 2 further elaborates on these details.

Our main theorem establishes that up to isomorphism, 2-dimensional endo-commutative straight algebras of type I over any
field $K$ are classified into a family $\{S(p, 0, 0, 0, 0, 0) \}_{ p\in\mathcal H}$, each of which is the structure matrix determined 
by the multiplication tables
\[
\begin{pmatrix}f&0\\0&pe\end{pmatrix}\, \, (p\in\mathcal H)
\]
with respect to a linear base $\{e, f\}$. Here, $\mathcal H$ denotes a complete representative system of $K^*/\approx$ for some equivalence relation $\approx$ of the multiplicative group $K^*=K\backslash\{0\}$. The details can be found in Theorem 1.  Furthermore, as applications of the main theorem, we show that if $K$ is cube-rootable, that is, $K=\{k^3 : k\in K\}$, then any 2-dimensional endo-commutative straight algebra of type I is isomorphic to the algebra $S(1, 0, 0, 0, 0, 0)$.  
Additionally, we establish that the set of all 2-dimensional endo-commutative straight algebras of type I over the rational number field is countably infinite, up to isomorphism.
A comprehensive description of these applications is provided in Corollaries 3 and 4.
Finally, we explore the isomorphism between type I and the other types.
The details are described in the final section.
The remaining challenges involve complete classifications of the algebras of type II and III.
Although these problems are complex, we anticipate that they hold simple and elegant solutions.

Throughout this paper, $K$ is an arbitrary field. 
In \cite{TST2,TST3}, it is assumed that $K$ is a non-trivial field, but since this paper focuses on the study of straight algebras, $K$ can be considered
as any field.

\section{Isomorphism criterion for 2-dimensional algebras and characterization of 2-dimensional endo-commutative algebras}\label{sec:iso-criterion}

For any $X=\small{\begin{pmatrix}a&b\\c&d\end{pmatrix}}\in GL_2(K)$, define
\[
\widetilde{X}=\begin{pmatrix}a^2&b^2&ab&ab\\c^2&d^2&cd&cd\\ac&bd&ad&bc\\ac&bd&bc&ad\end{pmatrix},
\]
where $GL_n(K)$ is the general linear group consisting of nonsingular $n\times n$ matrices over $K$.  Then the following result can be found in \cite{TST2}.

\begin{lem}
The mapping $X\mapsto\widetilde{X}$ is a group-isomorphism from $GL_2(K)$ into $GL_4(K)$ with $|\widetilde X|=|X|^4$.
\end{lem}

Let $A$ be a 2-dimensional algebra over $K$ with a linear base $\{e, f\}$.  Then $A$ is determined by the multiplication table $\begin{pmatrix}e^2&ef\\fe&f^2\end{pmatrix}$.  We write
\[
\left\{
    \begin{array}{@{\,}lll}
     e^2=a_1e+b_1f\\
     f^2=a_2e+b_2f\\
     ef=a_3e+b_3f\\
     fe=a_4e+b_4f\\
   \end{array}
  \right. 
\]
with $a_i, b_i\in K\, \, (1\le i\le4)$ and the matrix
$\begin{pmatrix}a_1&b_1\\a_2&b_2\\a_3&b_3\\a_4&b_4\end{pmatrix}$ is called the structure matrix of $A$ with respect to the base $\{e, f\}$. 

We hereafter will freely use the same symbol $A$ for the matrix and for the algebra because the algebra $A$ is determined by its structure matrix.  Then the following result can be found in \cite{TST2}, and this shows the isomorphism criterion for 2-dimensional algebras.

\begin{pro}
Let $A$ and $A'$ be 2-dimensional algebras over $K$.  Then $A$ and $A'$ are isomorphic iff there is $X\in GL_2(K)$ such that 
\begin{equation}
A'=\widetilde{X^{-1}}AX.
\end{equation}
\end{pro}

\begin{cor} Let $A$ and $A'$ be 2-dimensional algebras over $K$.  If $A$ and $A'$ are isomorphic, then ${\rm{rank}}\, A={\rm{rank}}\, A'$.
\end{cor}

When (1) holds, we say that the matrices $A$ and $A'$ are equivalent and refer to $X$ as a $\it{transformation\, \, matrix}$ for the equivalence $A\cong A'$.  Also, we call this $X$ a transformation matrix for the isomorphism $A\cong A'$ as well.
\vspace{2mm}

For any $a_1, b_1, a_2, b_2, a_3, b_3, a_4, b_4\in K$, we consider the following system consisting of eight cubic equations:
\begin{equation}
\left\{\begin{array}{@{\,}lll}   
a_1^2a_2+b_1a_2b_2+a_1b_2a_3+b_1a_2a_4=a_1a_3^2+a_2b_3^2+a_3^2b_3+a_3b_3a_4\\
a_1^2a_2+b_1a_2b_2+b_1a_2a_3+a_1b_2a_4=a_1a_4^2+a_2b_4^2+a_3a_4b_4+a_4^2b_4\\
a_1^2a_4+b_1a_4^2+b_1a_2b_4+a_1a_3b_4=a_1^2a_3+b_1a_2b_3+b_1a_3^2+a_1b_3a_4\\
a_2(a_1a_4+a_4b_4+b_2b_4)=a_2(a_1a_3+b_2b_3+a_3b_3)\\
a_1b_1a_2+b_1b_2^2+a_1b_2b_3+b_1a_2b_4=b_1a_3^2+b_2b_3^2+a_3b_3^2+a_3b_3b_4\\
a_1b_1a_2+b_1b_2^2+b_1a_2b_3+a_1b_2b_4= b_1a_4^2+b_2b_4^2+b_3a_4b_4+a_4b_4^2\\
b_1(a_1a_4+a_4b_4+b_2b_4)=b_1(a_1a_3+b_2b_3+a_3b_3)\\
b_1a_2a_4+b_2b_3a_4+b_2^2b_4+a_2b_4^2=b_1a_2a_3+b_2^2b_3+a_2b_3^2+b_2a_3b_4.  
\end{array} \right. 
\end{equation}
The following result can be found in \cite{TST2}, and this gives a necessary and sufficient condition for 2-dimensional algebra $A$ over $K$ to be endo-commutative.

\begin{pro}
Let $A$ be a 2-dimensional algebra over $K$ with structure matrix $\small{\begin{pmatrix}a_1&b_1\\a_2&b_2\\a_3&b_3\\a_4&b_4\end{pmatrix}}$.  Then $A$ is endo-commutative iff  the scalars $a_1, b_1, a_2, b_2, a_3, b_3, a_4, b_4$ satisfy the system $\rm{(2)}$ consisting of eight cubic equations.
\end{pro}

\section{2-dimensional endo-commutative straight algebras}\label{straight}

Recall that a 2-dimensional algebra is curled if the square of any element $x$ is a scalar multiplication of $x$, otherwise, it is straight.  For a point $(p, q, a, b, c, d)\in K^6$, let $S(p, q, a, b, c, d)$ be an algebra over $K$ with the multiplication table
\[
\begin{pmatrix}f&ae+bf\\ce+df&pe+qf\end{pmatrix}
\]
with respect to the linear base $\{e, f\}$.  The structure matrix of this algebra is given by 
\[
S(p, q, a, b, c, d)=\begin{pmatrix}0&1\\p&q\\a&b\\c&d\end{pmatrix}.
\]

Of course, this algebra is straight.  Conversely, we see easily  that any 2-dimensional straight algebra over $K$ is isomorphic to some $S(p, q, a, b, c, d)$ by replacing the bases.  The following results can be found in \cite{TST3}.

\begin{lem}
The straight algebra $S(p, q, a, b, c, d)$ is endo-commutative iff the point $(p, q, a, b, c, d)\in K^6$ satisfies 
\begin{equation}
\left\{\begin{array}{@{\,}lll}   
pq+pc=pb^2+a^2b+abc\\
p(c-a)=(b-d)\{p(b+d)-q(a+c)\}\\
p(d-b)=a^2-c^2\\
q^2+pd=a^2+qb^2+ab^2+abd\\
q(d-b)=ab-cd.\\  
\end{array} \right. 
\end{equation}
\end{lem}

\begin{lem}
The straight algebras $S(p, q, a, b, c, d)$ and $S(p', q', a', b', c', d')$ are isomorphic iff there are $x, y, z, w\in K$ with $\left|\begin{matrix}x&y\\z&w\end{matrix}\right|\ne0$ such that 
\begin{equation}
\left\{\begin{array}{@{\,}lll}
p'y^2+(a'+c')xy=z\\
x^2+q'y^2+(b'+d')xy=w\\
p'w^2+(a'+c')zw=px+qz\\
z^2+q'w^2+(b'+d')zw=py+qw\\
p'yw+a'xw+c'yz=ax+bz\\
xz+q'yw+b'xw+d'yz=ay+bw\\
p'yw+a'yz+c'xw=cx+dz\\
xz+q'yw+b'yz+d'xw=cy+dw
\end{array} \right. 
\end{equation}
holds.
\end{lem}

In Lemma 3, note that $X:=\begin{pmatrix}x&y\\z&w\end{pmatrix}$ is a transformation matrix for the isomorphism $S(p, q, a, b, c, d)\cong S(p', q', a', b', c', d')$.

\section{2-dimensional endo-commutative straight algebras of rank 2}\label{rank2}

Let $\mathcal{EC}$ be the family of all 2-dimensional endo-commutative algebras over $K$ with a linear base $\{e, f\}$.  Define 
\[
\mathcal{ECS}_{{\rm{(rank \, 2)}}}=\{S(p, q, a, b, c, d)\in\mathcal{EC} : {\rm{rank}}\, S(p, q, a, b, c, d)=2\}.
\]
Here $S(p, q, a, b, c, d)\in\mathcal{EC}$ means that the algebra $S(p, q, a, b, c, d)$ is endo-commutative.
Obviously, the rank of $S(p,q,a,b,c,d)$ is 2 iff it is not the case that all of $p$, $a$, and $c$ are 0. Therefore, we categorize 
$\mathcal{ECS}_{{\rm{(rank \, 2)}}}$ into three types based on the number of nonzero entries of $p$, $a$, and $c$.

Define
\[
\mathcal{ECS}_{{\rm{I}}}=\mathcal{ECS}_{001}\cup\mathcal{ECS}_{010}\cup\mathcal{ECS}_{100},
\]
where 
\[
\left\{\begin{array}{@{\,}lll}  
\mathcal{ECS}_{001}=\{S(p, q, a, b, c, d)\in\mathcal{EC} : p=a=0, c\ne0\}\\
\mathcal{ECS}_{010}=\{S(p, q, a, b, c, d)\in\mathcal{EC} : p=c=0, a\ne0\}\\
\mathcal{ECS}_{100}=\{S(p, q, a, b, c, d)\in\mathcal{EC} : a=c=0, p\ne0\}.
\end{array} \right. 
\]
Then $\mathcal{ECS}_{{\rm{I}}}$ is a subfamily of $\mathcal{ECS}_{{\rm{(rank \, 2)}}}$.  We say that each algebra in $\mathcal{ECS}_{{\rm{I}}}$ is of type I.   Define
\[
\mathcal{ECS}_{{\rm{II}}}=\mathcal{ECS}_{011}\cup\mathcal{ECS}_{101}\cup\mathcal{ECS}_{110},
\]
where 
\[
\left\{\begin{array}{@{\,}lll}  
\mathcal{ECS}_{011}=\{S(p, q, a, b, c, d)\in\mathcal{EC} : p=0, a, c\ne0\}\\
\mathcal{ECS}_{101}=\{S(p, q, a, b, c, d)\in\mathcal{EC} : a=0, p, c\ne0\}\\
\mathcal{ECS}_{110}=\{S(p, q, a, b, c, d)\in\mathcal{EC} : c=0, p, a\ne0\}.
\end{array} \right. 
\]
Then $\mathcal{ECS}_{{\rm{II}}}$ is a subfamily of $\mathcal{ECS}_{{\rm{(rank \, 2)}}}$.  We say that each algebra in $\mathcal{ECS}_{{\rm{II}}}$ is of type II.  Define
\[
\mathcal{ECS}_{{\rm{III}}}=\{S(p, q, a, b, c, d)\in\mathcal{EC} : p, a, c\ne0\}.
\]
Then $\mathcal{ECS}_{{\rm{III}}}$ is a subfamily of $\mathcal{ECS}_{{\rm{(rank \, 2)}}}$.  We say that each algebra in $\mathcal{ECS}_{{\rm{III}}}$ is of type III.
Obviously, we have the following disjoint union:
\[
\mathcal{ECS}_{{\rm{(rank \, 2)}}}=\mathcal{ECS}_{{\rm{I}}}\sqcup\mathcal{ECS}_{{\rm{II}}}\sqcup\mathcal{ECS}_{{\rm{III}}}.
\]

\section{Classification of $\mathcal{ECS}_{{\rm{I}}}$}\label{typeI}

The family of  2-dimensional endo-commutative straight algebras of type I is characterized as follows:
\begin{pro}
$\mathcal{ECS}_{{\rm{I}}}=\{S(p, 0, 0, 0, 0, 0) : p\ne0\}$.
\end{pro}

\begin{proof}
(i) $\mathcal{ECS}_{001}=\emptyset$.  Suppose $\mathcal{ECS}_{001}\ne\emptyset$, and take $S(p, q, a, b, c, d)\in\mathcal{ECS}_{001}$.  Since $S(p, q, a, b, c, d)$ is endo-commutative, the point $(p, q, a, b, c, d)\in K^6$ satisfies (3) by Lemma 2.  Also since $p=a=0$ and $c\ne0$, we see that (3) is rewritten as
\[
\left\{\begin{array}{@{\,}lll}   
0=-qc(b-d)\\
0=-c^2\\
q^2=qb^2\\
q(d-b)=-cd.\\  
\end{array} \right.
\]
However, the second equation contradicts $c\ne0$, hence the assertion (i) holds.
\vspace{2mm}

(ii) $\mathcal{ECS}_{010}=\emptyset$.  Suppose $\mathcal{ECS}_{010}\ne\emptyset$, and take $S(p, q, a, b, c, d)\in\mathcal{ECS}_{010}$.  Since $S(p, q, a, b, c, d)$ is endo-commutative, the point $(p, q, a, b, c, d)\in K^6$ satisfies (3) by Lemma 2.  Also since $p=c=0$ and $a\ne0$, we see that (3) is rewritten as
\[
\left\{\begin{array}{@{\,}lll}   
0=a^2b\\
0=-qa(b-d)\\
0=a^2\\
q^2=a^2+qb^2+ab^2+abd\\
q(d-b)=ab. 
\end{array} \right. 
\]
However, the third equation contradicts $a\ne0$, hence the assertion (ii) holds.
\vspace{2mm}

(iii) $\mathcal{ECS}_{100}=\{S(p, 0, 0, 0, 0, 0) : p\ne0\}$.  Any point $(p, 0, 0, 0, 0, 0)$ in $K^6$ with $p\ne0$ clearly satisfies (3), and hence $S(p, 0, 0, 0, 0, 0)$ is endo-commutative by Lemma 2, so $\{S(p, 0, 0, 0, 0, 0) : p\ne0\}\subseteq\mathcal{ECS}_{100}$.  

To show the opposite inclusion, take $S(p, q, a, b, c, d)\in\mathcal{ECS}_{100}$ arbitrarily.  Since $S(p, q, a, b, c, d)$ is endo-commutative, the point $(p, q, a, b, c, d)\in K^6$ satisfies (3) by Lemma 2 again.  Also since $a=c=0$ and $p\ne0$, we see that (3) is rewritten as
\[
\left\{\begin{array}{@{\,}lll}   
pq=pb^2\\
0=(b-d)p(b+d)\\
p(d-b)=0\\
q^2+pd=qb^2\\
q(d-b)=0,  
\end{array} \right. 
\]
which is rewritten as $\left\{\begin{array}{@{\,}lll}   
q=b^2\\
b=d\\
q^2+pb=qb^2,
\end{array} \right.$ hence $q=b=d=0$ by a simple calculation.  Therefore  we see that $(p, q, a, b, c, d)=(p, 0, 0, 0, 0, 0)$ and $p\ne0$.  Thus we see that the assertion (iii) holds.  

Since $\mathcal{ECS}_{{\rm{I}}}=\mathcal{ECS}_{001}\cup\mathcal{ECS}_{010}\cup\mathcal{ECS}_{100}$, the proposition follows from (i), (ii) and (iii) above.
\end{proof}

\begin{cor}
Any 2-dimensional endo-commutative straight algebra of type I over $K$ is commutative, but neither unital nor associative. 
\end{cor}

\begin{proof}
  Let $A$ be a 2-dimensional endo-commutative straight algebra of type I over $K$.  By Proposition 3, there is $p\in K$ with $p\ne0$ such that $A$ has a multiplication table $\begin{pmatrix}f&0\\0&pe\end{pmatrix}$ with respect to a linear base $\{e, f\}$.  Since this table is symmetric, it is clear that $A$ is commutative.
    If $A$ has an identity element $u=\alpha e+\beta f$, where $\alpha, \beta\in K$, then we have $e=ue=\alpha e^2+\beta fe=\alpha f$, a contradiction, so $A$ is non-unital.  Also since $e^2f=f^2=pe\ne0$ and $e(ef)=e0=0$, it follows that $e^2f\ne e(ef)$, hence $A$ is non-associative.
\end{proof}

We define
\[
(K^*)^3=\{k^3 : k\in K^*\},
\]
where $K^*:=K\, \backslash\, \{0\}$ is the multiplicative group of $K$.  
Then $(K^*)^3$ forms a subgroup of $K^*$.  We define an equivalence relation $\sim$ in $K^*$
as follows: for any $a, b\in K^*$, $a\sim b$ iff $a/b\in (K^*)^3$.  Moreover, for $a, b\in K^*$, write $a\approx b$ if $a\sim b$ or $a^2\sim b$.  Then,
$\approx$ is also an equivalence relation in $K^*$.  In fact, it is obvious that for any $a\in (K^*)^3$, $a\approx a$. If $a\approx b$, then $a\sim b$ or $a^2\sim b$.
When $a\sim b$, we have $b\sim a$ and so $b\approx a$. When $a^2\sim b$, we have $b\sim a^2$, and hence $b^2\sim a^4\sim a$ which implies $b\approx a$.
Suppose that the scalars $a, b, c\in K^*$ satisfy $a\approx b$ and $b\approx c$.  Then [$a\sim b$ or $a^2\sim b$] and [$b\sim c$ or $b^2\sim c$] hold.

(i) In the case where $a\sim b$ and $b\sim c$, it is clear that $a\sim c$.

(ii) In the case where $a\sim b$ and $b^2\sim c$, since $a^2\sim b^2$, it follows that $a^2\sim c$.

(iii) In the case where $a^2\sim b$ and $b\sim c$, it is clear that $a^2\sim c$.

(iv) In the case where $a^2\sim b$ and $b^2\sim c$, there are $x, y\in K^*$ such that $a^2=bx^3$ and $b^2=cy^3$. Then $(a^2/x^3)^2=cy^3$, and hence $a^4=cx^6y^3$.  
Hence we have $a=cx^6y^3/a^3=cz^3$, where $z=x^2y/a$, that is, $a\sim c$.  By (i), (ii) (iii) and (iv), we see $a\approx c$.

\begin{lem}
 Let $p, p'\in K^*$.  Then $S(p, 0, 0, 0, 0, 0)\cong S(p', 0, 0, 0, 0, 0)$ iff $p\approx p'$.
\end{lem}

\begin{proof}
Substituting $q=a=b=c=d=0$ and $q'=a'=b'=c'=d'=0$ into (4), we obtain
\begin{equation}
\left\{\begin{array}{@{\,}lll}
p'y^2=z\\
x^2=w\\
p'w^2=px\\
z^2=py\\
yw=0\\
xz=0.
\end{array} \right. 
\end{equation}
Therefore, by Lemma 3, $S(p, 0, 0, 0, 0, 0)\cong S(p', 0, 0, 0, 0, 0)$ iff there are $x, y, z, w\in K$ with $\left|\begin{matrix}x&y\\z&w\end{matrix}\right|\ne0$ such that  (5) holds.  
\vspace{2mm}

(i) Suppose $p\sim p'$.  Then there is $x\in K^*$ with $p=p'x^3$.  Put $y=z=0$ and $w=x^2$.  Then $\left|\begin{matrix}x&y\\z&w\end{matrix}\right|=x^3\ne0$.  Also since $p'w^2=p'x^4=p'x^3x=px$, we see that the scalars $x, y, z, w$ satisfy (5), and hence $S(p, 0, 0, 0, 0, 0)\cong S(p', 0, 0, 0, 0, 0)$ with a transformation matrix $\begin{pmatrix}x&y\\z&w\end{pmatrix}$.  
\vspace{2mm}

(ii) Suppose $p^2\sim p'$.  Then there is $z\in K^*$ with $p^2=p'z^3$.  Put $x=w=0$ and $y=z^2/p$. Then $\left|\begin{matrix}x&y\\z&w\end{matrix}\right|=-yz=-z^3/p\ne0$. Also since $p'y^2=p'z^4/p^2=p'z^3z/p^2=p^2z/p^2=z$ and $z^2=py$, we see that the scalars $x, y, z, w$ satisfy (5), and hence $S(p, 0, 0, 0, 0, 0)\cong S(p', 0, 0, 0, 0, 0)$ with a transformation matrix $\begin{pmatrix}x&y\\z&w\end{pmatrix}$. 
\vspace{2mm}

(iii) Suppose $S(p, 0, 0, 0, 0, 0)\cong S(p', 0, 0, 0, 0, 0)$.  Then there are $x, y, z, w\in K$ with $\left|\begin{matrix}x&y\\z&w\end{matrix}\right|\ne0$ such that  (5) holds.  

(iii-1) The case where $x\ne0$.  By the second and third equations of (5), we obtain $p'x^4=px$, hence $p'x^3=p$, that is, $p\sim p'$, so $p\approx p'$.

(iii-2) The case where $x=0$.  In this case, we see $y\ne0$ because $\left|\begin{matrix}x&y\\z&w\end{matrix}\right|=-yz\ne0$.  By the first and fourth equations of (5), we obtain $p'^2y^4=py$, hence $p'^2y^3=p$, that is, $p^2\sim p'$, so $p\approx p'$.

By (i), (ii), (iii), we obtain the desired result.
\end{proof}

Take $\mathcal H\subseteq K^*$ as a complete representative system of the equivalent classes $K^*/\approx$.
From Proposition 3 and Lemma 4, we obtain the following result.

\begin{thm}
Up to isomorphism, 2-dimensional endo-commutative straight algebras of type I over $K$ are classified into the family $\{S(p, 0, 0, 0, 0, 0) \}_{ p\in\mathcal H}$ defined by the multiplication tables
\[
\begin{pmatrix}f&0\\0&pe\end{pmatrix}\, \, (p\in\mathcal H)
\]
with respect to a linear base $\{e, f\}$.
\end{thm}

We refer to $K$ as cube-rootable if $K=\{k^3 : k\in K\}$.  For example, the complex number field $\mathbb C$, the real number field $\mathbb R$, $\mathbb Z_2$ and $\mathbb Z_3$ are cube-rootable.  However, the rational number field $\mathbb Q$ is not cube-rootable.

As an application of Theorem 1, we have the following.

\begin{cor}
Suppose that $K$ is cube-rootable.  Then all 2-dimensional endo-commutative straight algebras of type I over $K$ are isomorphic to the algebra\\ $S(1, 0, 0, 0, 0, 0)$ with the multiplication table $\begin{pmatrix}f&0\\0&e\end{pmatrix}$.
\end{cor}
\vspace{1mm}

Corollary 3 asserts that when $K$ is cube-rootable, there exists a unique 2-dimensional endo-commutative straight algebra of type I over $K$, up to isomorphism.
However, it should be noted that depending on the choice of $K$, there can be infinitely many such algebras. First, we need the following lemma.
\vspace{1mm}

\begin{lem} $\sharp ({\mathbb Q}^*/\approx)=\infty$.
  \end{lem}
  
\begin{proof} We will show that every prime number can be taken as a representative of an equivalence class in ${\mathbb Q}^*/\approx$, ensuring that
  distinct primes correspond to distinct representatives. This completes the proof since the set of all primes is infinite.
  Therefore, our goal is to prove that if $p$ and $q$ are primes with $p\neq q$, then $p\not\approx q$, that is, $p\not\sim q$ and $p^2\not\sim q$.

  \noindent  (i) $p\not\sim q$.

  Assume $p\sim q$. Then, there exist integers $m$ and $n$ such that \[\dfrac{p}{q}=\left(\dfrac{n}{m}\right)^3.\]
  Thus, $pm^3=qn^3$. 
  Suppose that the prime factorizations of $m$ and $n$ are $p_1^{e_1}p_2^{e_2}\cdots p_r^{e_r}$ and $q_1^{f_1}q_2^{f_2}\cdots q_s^{f_s}$, respectively.
  Hence, \[pp_1^{3e_1}p_2^{3e_2}\cdots p_r^{3e_r}=qq_1^{3f_1}q_2^{3f_2}\cdots q_s^{3f_s}.\]
  Therefore, $p$ divides the quantity on the right hand side.
  Since $p\neq q$, $p$ must be either one of $q_1, q_2,\dots,q_s$, say $p=q_1$.
  By the uniqueness of factorization, the exponent of $p$ on the left hand side equals to the exponent of $p$ on the right hand side.
  If one of $p_1,p_2,\dots,p_r$ coincides with $p$, say $p_1=p$, then the exponent of $p$ on the left hand side is $3e_1+1$.
  On the other hand, the exponent of $p$ on the right hand side is $3f_1$, and so $3e_1+1=3f_1$, which is a contradiction.
  If none of $p_1,p_2,\dots,p_r$ coincides with $p$, by comparing the exponents of $p$ on both sides, we obtain $1=3f_1$, which is also a contradiction.

  \noindent  (ii) $p^2\not\sim q$.

  Assume $p^2\sim q$. Then, there exist integers $m$ and $n$ such that \[\dfrac{p^2}{q}=\left(\dfrac{n}{m}\right)^3.\]
  Thus, $p^2m^3=qn^3$. 
  Suppose that the prime factorizations of $m$ and $n$ are $p_1^{e_1}p_2^{e_2}\cdots p_r^{e_r}$ and $q_1^{f_1}q_2^{f_2}\cdots q_s^{f_s}$, respectively.
  Hence, \[p^2p_1^{3e_1}p_2^{3e_2}\cdots p_r^{3e_r}=qq_1^{3f_1}q_2^{3f_2}\cdots q_s^{3f_s}.\] Since $p\neq q$, $p$ must be either one of $q_1, q_2,\dots,q_s$, say $p=q_1$.
  By applying a similar argument as above, we can deduce that either $3e_1+2=3f_1$ or $2=3f_1$, both of which lead to a contradiction.
\end{proof}

\begin{cor}
The set of all 2-dimensional endo-commutative straight algebras of type I over $\mathbb Q$ is countably infinite up to isomorphism.
\end{cor}

\begin{proof}
As observed in the proof of Lemma 5, it follows that if $p$ and $q$ are primes
  with $p\neq q$, then $p\not\approx q$, and hence $S(p, 0, 0, 0, 0, 0)\ncong S(q, 0, 0, 0, 0, 0, )$ by Theorem 1.  This completes the proof.
\end{proof}

\section{Relationship of type I to types II and III}\label{relation}
\vspace{2mm}

For $n\in\mathbf N$, we refer to $K$ as $n$-rootable if $K=\{k^n : k\in K\}$.  Of course, \lq\lq cube-rootable" and \lq\lq 3-rootable" have the same meaning.  Also \lq\lq 2-rootable" is usually called \lq\lq square-rootable" or \lq\lq  quadratically closed".  Moreover, $K$ is said to be quadratically non-extendable if any quadratic polynomial in $K[x]$ has a root in $K$. If ${\rm{char}}\, K=2$, then $K$ is quadratically non-extendable iff $K$ is quadratically closed, but otherwise the \lq\lq if'' part is not necessarily true  (cf. \cite{SYTT1}). 

\begin{lem}
Let $S(p, q, a, b, c, d)\in\mathcal{ECS}_{\rm{II}}\cup\mathcal{ECS}_{\rm{III}}$ and $S(p', q', a', b', c', d')\in \mathcal{ECS}_{\rm{I}}$.  If $S(p, q, a, b, c, d)\cong S(p', q', a', b', c', d')$, then 
\[
(\sharp_1)\left\{\begin{array}{@{\,}lll}
p=0\\
q=-a\\
a\ne0\\
b=0\\
c=a\\
d=0
\end{array} \right.\, \, {\rm{or}}\, \,\, \,  (\sharp_2)\left\{\begin{array}{@{\,}lll}
p\ne0\\
q=-a\\
a\ne0\\
b=0\\
c=a\\
d=0.
\end{array} \right. 
\]
\end{lem}

\begin{proof}

Since $S(p, q, a, b, c, d)$ is of type II or III, we see from definition that at least two of $p, a$ and $c$ are nonzero.  Also since  $S(p', q', a', b', c', d')$ is of type I, it follows from Proposition 3 that $p'\ne0, q'=a'= b'=c'=d'=0$. 

Suppose $S(p, q, a, b, c, d)\cong S(p', q', a', b', c', d')$.  By Lemma 3, there are $x, y, z, w\in K$ with $\left|\begin{matrix}x&y\\z&w\end{matrix}\right|\ne0$ such that (4) holds.  Because $q'=a'= b'=c'=d'=0$, (4) is rewritten as
\begin{equation}
\left\{\begin{array}{@{\,}lll}
p'y^2=z\cdots\mbox{(6$-$1)}\\
x^2=w\cdots\mbox{(6$-$2)}\\
p'w^2=px+qz\cdots\mbox{(6$-$3)}\\
z^2=py+qw\cdots\mbox{(6$-$4)}\\
p'yw=ax+bz\cdots\mbox{(6$-$5)}\\
xz=ay+bw\cdots\mbox{(6$-$6)}\\
p'yw=cx+dz\cdots\mbox{(6$-$7)}\\
xz=cy+dw\cdots\mbox{(6$-$8)}.
\end{array} \right. 
\end{equation}
Then we shall show that $x, y, z, w\ne0$ and $b=d=0$ in (a)$\sim$(e) as follows:

(a) $w\ne0$.  Suppose $w=0$.  By (6$-$2), $x=0$.  Then $0\ne\left|\begin{matrix}x&y\\z&w\end{matrix}\right|=-yz$, hence $y\ne0, z\ne0$.   Therefore, $a=0$ by (6$-$6) and $c=0$ by (6$-$8), a contradiction.  
\vspace{2mm}

(b) $x\ne0$.  This follows from (6$-$2) and (a).
\vspace{2mm}

(c) $z\ne0$.  Suppose $z=0$.  By (6$-$1), $y=0$ because $p'\ne0$.  Then by (6$-$5), $0=ax$, hence $a=0$ because $x\ne0$.  Similarly, by (6$-$7), $0=cx$, hence $c=0$.  Then we obtain $a=c=0$, a contradiction.
\vspace{2mm}

(d) $y\ne0$.  This follows from (6$-$1) and (c).
\vspace{2mm}

(e) $b=d=0$.  Computing (6$-$5)$\times y-$ (6$-$6)$\times x$, we get  $p'y^2w-x^2z=b(zy-xw)$.  Also, computing (6$-$7)$\times y-$ (6$-$8)$\times x$, we get $p'y^2w-x^2z=d(zy-wx)$.  But we see $p'y^2w-x^2z=p'y^2x^2-x^2p'y^2=0$ by (6$-$1) and (6$-$2).  Then we obtain $b=d=0$ because $\left|\begin{matrix}x&y\\z&w\end{matrix}\right|\ne0$.
\vspace{2mm}

Substituting  (6$-$1) $z=p'y^2$ and (6$-$2) $w=x^2$ into (6$-$3)$\sim$(6$-$8), (6) is rewritten as 
\begin{equation}
\left\{\begin{array}{@{\,}lll}
p'y^2=z\cdots\mbox{(7$-$1)}\\
x^2=w\cdots\mbox{(7$-$2)}\\
p'x^4=px+qp'y^2\cdots\mbox{(7$-$3)}\\
p'^2y^4=py+qx^2\cdots\mbox{(7$-$4)}\\
p'yx=a\cdots\mbox{(7$-$5)}\\
p'yx=c\cdots\mbox{(7$-$6)}\\
\end{array} \right. 
\end{equation}
because $b=d=0$.  Substituting (7$-$5) $y=a/(p'x)$ into (7$-$3) and (7$-$4), (7) is rewritten as  
\begin{equation}
\left\{\begin{array}{@{\,}lll}
p'y^2=z\cdots\mbox{(8$-$1)}\\
x^2=w\cdots\mbox{(8$-$2)}\\
p'^2x^6=pp'x^3+qa^2\cdots\mbox{(8$-$3)}\\
a^4=pp'ax^3+qp'^2x^6\cdots\mbox{(8$-$4)}\\
p'xy=a\cdots\mbox{(8$-$5)}\\
a=c.\cdots\mbox{(8$-$6)}\\
\end{array} \right. 
\end{equation}
We note that $a\ne0$ by (8$-$5).  Since $S(p, q, a, b, c, d)$ is endo-commutative, it follows from Lemma 2 that (3) holds.  Since $a=c$ and $b=d=0$ by (8$-$6) and (e), (3) is rewritten as
\begin{equation}
\left\{\begin{array}{@{\,}lll}   
p(q+a)=0\\
q^2=a^2,
\end{array} \right. 
\end{equation}
which is equivalent to $\left\{\begin{array}{@{\,}lll}   
p=0\\
a=q
\end{array} \right. $ or $\left\{\begin{array}{@{\,}lll}   
p=0\\
a=-q
\end{array} \right. $or $\left\{\begin{array}{@{\,}lll}   
p\ne0\\
q=-a.
\end{array} \right. $  However, the case  $\left\{\begin{array}{@{\,}lll}   
p=0\\
a=q
\end{array} \right. $ does not occur.  In fact, if it occurs, we see from (8$-$3) that $p'^2x^6=a^3$, hence 
\[
\left|\begin{matrix}x&y\\z&w\end{matrix}\right|=xw-yz=x^3-\frac{a}{p'x}\, p'\left(\frac{a}{p'x}\right)^2=x^3-\frac{a^3}{p'^2x^3}=x^3-x^3=0,
\]
a contradiction.  Thus we get either ${\rm{(\sharp_1)}}$ or ${\rm{(\sharp_2)}}$. 
\end{proof}

  We can derive the following statements from Lemma 6, but we will omit the proofs as they deviate from the main theme of this paper.
\vspace{2mm}

(I) If ${\rm{char}}\, K=2$, then any algebra of type I is not isomorphic to any algebra of type II.
\vspace{2mm}

(II) If ${\rm{char}}\, K\ne2$ and $K$ is square-rootable, then there exists an algebra of type II that is isomorphic to an algebra of type I.
\vspace{2mm}

(III) If ${\rm{char}}\, K=2$ and $K$ is both quadratically non-extendable and cube-rootable, then there exists an algebra of type III
that is isomorphic to an algebra of type I.
\vspace{2mm}

(IV) If ${\rm{char}}\, K\ne2$ and $K$ is 6-rootable, then there exists an algebra of type III that is isomorphic to an algebra of type I.
\vspace{2mm}

Initially, we believed that algebras of type I, type II, and type III are not isomorphic to each other. However, the above-mentioned facts have proven the assumption
wrong. Consequently, it is expected that the classifications of type II and type III algebras will pose significant challenges in future studies.

\end{document}